\documentclass[11pt]{article}
\usepackage{graphicx}
\usepackage{subfigure}
\usepackage{multirow}
\usepackage{url}
\usepackage{amssymb}
\usepackage{amsmath}
\usepackage{amsfonts}
\usepackage{wrapfig}
\usepackage{afterpage}
\usepackage{float}
\usepackage{algorithm}
\usepackage[noend]{algpseudocode}
\usepackage{epsfig}
\usepackage{epstopdf}
\usepackage{longtable}
\usepackage{url} 
\usepackage{caption}
\usepackage{listings}
\usepackage{color} 
\usepackage{makecell}
\definecolor{mygreen}{RGB}{28,172,0} 
\definecolor{mylilas}{RGB}{170,55,241}
\newcommand{\real}{\mathbb{R}}
\usepackage{tabularx}
\usepackage{hyperref}

\usepackage{listings}
\usepackage{color} 
\definecolor{mygreen}{RGB}{28,172,0} 
\definecolor{mylilas}{RGB}{170,55,241}

\begin{document}
 \vspace{1cm}
\begin{center}
\Large {\bf  Topology Optimization with Bilevel Knapsack: An Efficient 51 Lines MATLAB Code }
 \vspace{0.5cm}\\
{\bf Vittorio Latorre}\footnote{Corresponding author. Email: {\em v.latorre@federation.edu.au}},   

{\em  Faculty of Science and Technology,\\ Federation University Australia, Mt Helen, Victoria 3353, Australia}\\
\end{center}

\begin{abstract}
This paper presents an efficient 51 lines Matlab code to solve topology optimization problems. By the fact that the presented code is based on an hard 0-1 optimization method that handles the integer part of the optimization in a simple fashion and in sub-polynomial time, it has been possible to not only reduce the number of lines to 51 without sacrificing any readability, but also to obtain void-solid designs with close to none checkerboard patterns with improved efficiency. 
The numerical results in the paper show that the proposed method has the best average times compared to several codes available in literature. 

We also present extensions to different boundary conditions and to three dimensional designs. The code can be used by students and the newcomers in topology optimization because of its simplicity and readability. The 51 lines Matlab code and the presented extensions can be downloaded on the webpage \url{https://github.com/vlatorre847/TOSSE}{}.
\end{abstract}
 
\section{Introduction}
Topology Optimization consists in the engineering problem of placing material into a structure considering the possible loads on such structure and the geometric and physical constrains in order to obtain the best structural performance. Such problem in mechanics was first introduced in \cite{bk88} that is considered the seminal paper for the subject, and several approaches has been developed in order to find its solution. 
Among the most popular methods we name the density approach \cite{b89}, that uses the concept of density instead of clear void-solid elements to reduce the complexity of the mathematical problem associated with the topology optimization. The density method known as SIMP (Simplified Isotropic Material with Penalization)  is among the most popular algorithms for topology optimization \cite{zr91,ml92,sig99,acs11} due to its efficiency. Similar to the density approach are the the phase field \cite{bc03} approach and the  topological derivatives approach \cite{sz01}. Other methods include the
  evolutionary   approach \cite{xs93,hy10} based on purely discrete strategies that start from completely filled structures and then decrease the percentage of volume that can be occupied,
hard-killing elements and  reintroducing them if considered rewarding, and the level set  approach \cite{aft02,wei18} based on the level set functions that define the contours of the topology,   

Recently in \cite{gao18} a new mathematical approach has been proposed based on the minimization of the total potential energy of the structure. In this paper the topology optimization problem has been formulated as a Bilevel Mixed Integer optimization problem, where at the higher level a knapsack optimization problem is solved in the mass density variables $\boldsymbol\rho$ and at the lower level a continuous optimization problem is solved in the displacement variable $\bold u$. The main difficulty of this formulation is that the knapsack problem is NP-Hard that is its global solution cannot be reached in polynomial time. 

In this paper we  use the formulation in \cite{gao18} and present an efficient 51 lines code written in Matlab for topology optimization. This code shares many similarities with the well know 88 lines code presented in \cite{acs11}, that can be considered as its starting point. However, the new approach is an hard 0-1 evolutionary optimization method, based on a different mathematical formulation, and with
 a better average  efficiency. This is possible not only because of the new formulation presented in \cite{gao18}, but also to a further assumption that makes the solution of the  knapsack problem trivial. As a matter of facts we assume that 
all the elements of the topology have the same size. In this way, the knapsack problem is greatly simplified and its solution can be found in sub-polynomial time. We underline that this assumption is often made in many educational codes that are the basis for more complex methods used in industrial softwares. The code we present is not only shorter than the 88 lines code, but also as quite as readable and easily comparable with other codes presented in literature, and therefore can be used for education purposes and for the practitioners who approach the topology optimization problem for the first time.
Furthermore the same size elements assumption is not too  restrictive, as there already exist methods capable of finding the solution of knapsack problems with tens of thousand of variables in a matter of milliseconds \cite{pis05}. 

Based on the bilevel formulation and the same size elements assumption, we present a Topology Optimization Same Size Elements (TOSSE) algorithm which main features consist in:
\begin{itemize}
\item  An hard $\{0,1\}$ optimization method, that returns solutions with clear void-solid elements without using any filter or any other artificial technique;
\item Almost no checkerboard patterns in the structures;
\item Good numerical efficiency.
\end{itemize}
We also present extensive numerical results on the Messerschmitt-B\"{o}lkow-Blohm (MBB) 
 beam with comparison with some of the most popular methodologies in the topology optimization literature, varying the size of the meshes of the finite elements method. We show that the TOSSE algorithm is able to generate topologies with clear void-solid elements and no checkerboard patterns similar to the ones yielded by the evolutionary methods in elapsed times comparable or smaller than the SIMP methods. 

This paper is associated with the github project: \url{https://github.com/vlatorre847/TOSSE}{}. On the web page is possible not only to download the 51 lines code, but also Matlab files containing the extension of this methodology to the 2D and 3D cantilever beams. 

The paper is organized as follows: in the next section we present the bilevel formulation and discuss it, then in Section \ref{mat} we  discuss the MATLAB implementation of the algorithm. In Section \ref{res} we present the numerical results and several examples, in Section \ref{ext} we present the extensions to  the two  and three dimensional  cantilever beams and in Section \ref{con} we report the conclusions.

\section{Problem Formulation}\label{prob_form}
Consider an elastically deformable body which reference domain is $\Omega\in\real^d$, $d=2,3$, with boundary $\Gamma=\partial \Omega$. This body is subjected to a body force $\textbf b$ in the reference domain $\Omega$ and a given surface traction $\bold{t}(\bold{x})$ of dead load type on a part of the boundary that we call $\Gamma_t\subset \Gamma$ while the body is fixed on the remaining $\Gamma_u=\Gamma\backslash \Gamma_t$ surface. The total potential energy of this body deformed by the body force $\textbf b$ and the surface traction $\bold{t}(\bold{x})$ is:
\begin{equation}\label{eq: intpon}
\Pi(\bold{u},\rho)=\int_{\Omega}W(\nabla \bold{u})\rho d\Omega-\int_{\Omega} \bold{u}\bold{b}\rho d\Omega-\int{\Gamma_t} \bold{u} \bold{t} d \Gamma,
\end{equation}
where $\bold{u}$ denotes the displacement, $\boldsymbol\rho$ denotes the mass density and it is a discrete design variable and $W(\nabla \bold{u})$ is the stored energy density of the deformation tensor $\nabla \bold{u}$.

By using an interpolation matrix $\bold{N}_e$, it is possible to use the finite element method to discretize the domain $\Omega$ according to a predefined mesh with $nelx$ elements on the x axis, $nely$ elements on the y axis and $nelz$ elements on the z axis for three dimensional structures. The discretized design $\{\Omega_e\}$ is composed by $n=nelx*nely$  elements for 2D structures and $n=nelx*nely*nelz$ for 3D structures, and the  discretized versions of the variables $\bold{u}$ and $\boldsymbol\rho$ in $\Omega_e$ are:
$$
\bold{u}(\bold{x})=\bold{N}_e(\bold{x})\bold{u}_e,\quad \boldsymbol\rho(\bold{x})=\rho_e \in \{ 0,1\} \forall \bold{x}\in \Omega_e,
$$
where $\bold{u}_e$ is the nodal displacement vector belonging to ${\cal U}_a\subset\real^m$, that is  the space of kinematically admissible displacement fields, where in general $m=2*(nelx+1)*(nely+1)$ for 2D structures and $m=3*(nelx+1)*(nely+1)*(nelz+1)$ for 3D structures, and the design variable $\boldsymbol\rho \in \{ 0,1\}^n$ determines if the $e^{th}$ element $\Omega_e$ is void, with $\rho_e=0$ or solid with $\rho_e=1$.

By using the finite element method it is possible to recast the (\ref{eq: intpon}) as a real valued function:
$$
\Pi_r(\bold{u},\boldsymbol\rho)=C(\bold{u},\boldsymbol\rho)-\bold{u}^T\bold{f}
$$
where 
$$
\bold{f}=\left\{\int_{\Gamma_t^e} \bold{N}_e(\bold{x})^T\bold{t}_e(\bold{x}) d\Gamma \right\}\in\real^m.
$$
and $C(\bold{u},\boldsymbol\rho)$ is the part of the potential energy that connects the displacement vector $\bold{u}$ to the design variables $\boldsymbol\rho$. 

It is well known that the variational analysis to find the deformation field is possible only if the design structure is given, while the structure can be determined if the deformation field is known. Therefore the variational analysis for the deformation field and the optimal structural design must be performed concurrently. Consequently the topology optimization problem for  linear elastic structures can be formulated as a bi-level integer nonlinear optimization problem:
\begin{equation}\label{eq: biknap}
\begin{array}{cc}
\displaystyle\min_{\boldsymbol\rho}&\displaystyle \Phi(\bold{u},\boldsymbol\rho)\\
\displaystyle s.t. & \displaystyle\sum_{e=1}^{n} v_e\rho_e\le V\vspace{0.5em}\\
&\displaystyle\bold{u}=\mbox{arg} \min_{\bold{v\in{\cal U}_a}}\left\{\Pi_r(\bold{v},\boldsymbol\rho)\right\}
\end{array}.
\end{equation}
In this formulation, function $\Phi(\bold{u},\boldsymbol\rho)$ represents the target for the upper level problem with  design variable $\boldsymbol\rho$, and  the problem consists in finding the optimal topological shape  under a knapsack constrain in which $v_e$ represents the size of the $e$-th element and $V$ is the desired volume. The lower level problem is in the displacement variable $\bold{u}$ and consists in assuring that the deformation is such that the lowest value of the potential energy for a given design is reached. 

The cost function $\Phi(\bold{u},\boldsymbol\rho)$ depends on the nature of the material and the method used to approach the design problem. As the upper level problem is in the design variable $\boldsymbol\rho$ with the displacement fields fixed by the displacement variables $\bold{u}$, the aim of this optimization process is to maintain the physical elements that have the highest potential energy stored in themselves. A way to express this relation is to set $\Phi(\bold{u},\boldsymbol\rho)=-\Pi_r(\bold{v},\boldsymbol\rho)$.

$\Phi(\bold{u},\boldsymbol\rho)$ should contain the numerical relations between the variables $\bold{u}$ and $\boldsymbol\rho$, or at least a reasonable approximation of such relations.
Therefore the main complexity of formulation (\ref{eq: biknap}) is the discrete $0-1$ nature of the design variables. In many approaches, like for example the popular SIMP approach \cite{bsbook}, $\boldsymbol\rho$ is considered continuous and some penalty is applied in order to steer the final solution toward discrete $0-1$ values. Our approach consists in an hard $0-1$ approach, however the difficulty remains in understanding how the design variables are influenced by the displacement variables and vice-versa. It is well know that the displacement fields are given implicitly in term of the design variable through the equilibrium equations, that is $\bold{u}$ must be a solution of the lower level problem, however the relation cannot be found explicitly given the discrete nature of $\boldsymbol\rho$. 

For this reason, we use the sensitivity analysis \cite{bsbook,hy10} to express this relation. The standard procedure in the sensitivity analysis is to consider the design problem, that is the higher level problem, as an optimization problem in the design variables only. The displacement variable $\bold{u}$ influences the design variables through the first order conditions of the total potential function $\Pi_r(\bold{u},\boldsymbol\rho)$  with respect to  $\boldsymbol\rho$.

In order to proceed we assume that the topology optimization is performed for a
 linear elastic structure without body force similarly to \cite{acs11,hy10} , with the total potential energy being a quadratic function:
\begin{equation}\label{eq: linelq}
\Pi_r(\bold{u},\boldsymbol\rho)=\frac{1}{2} \bold{u}^T \bold{K}(\boldsymbol\rho)\bold{u}-\bold{u}^T\bold{f}
\end{equation}
where $\bold{K}(\boldsymbol\rho)=\{\rho_e \bold{K}_e \}\in\real^{n\times n}$ is the overall stiffness matrix, obtained assembling the element-wise stiffness matrix $\rho_e \bold{K}_e$ for each element $\Omega_e$. 
%
For these materials the adjoint method for sensitivity analysis \cite{bsbook} gives the following derivative for  the $e^{th}$ element:
\begin{equation}\label{eq: costeq}
\frac{\partial\Pi_r(\bold{u},\boldsymbol\rho)}{\partial \rho_e}=\frac{1}{2}p\rho_e^{p-1} \bold{u}_e^T\bold{K}_e\bold{u}_e
\end{equation}
It is possible to notice the exponent $p$ for the design variable $\rho_e$. This exponent is generally used in topology optimization when $\boldsymbol\rho$ is considered continuous to operate the penalization that brings the solution to integer values.
Because of the hard $0-1$ formulation, the solution it is not affected for any $p>1$ in the (\ref{eq: costeq}).
Choosing $p=1$ would erase any dependence of the (\ref{eq: costeq}) from the design variable. However $\frac{\partial\Pi_r(\bold{u},\boldsymbol\rho)}{\partial \rho_e}$ is also effected by the values of the design variables and this effect is hidden in the displacement $\bold{u}$. Therefore we at least consider the influence of the element itself, set the value of $p=2$ and write:   
\begin{equation}\label{eq: c2}
\frac{\partial\Pi_r(\bold{u},\boldsymbol\rho)}{\partial \rho_e}=\rho_e \bold{u}_e^T\bold{K}_e\bold{u}_e.
\end{equation}
Consequently for linear elastic structures we use the derivative given in (\ref{eq: costeq}) to linearize the relation between $\bold{u}$ and $\boldsymbol\rho$   and write:
$$
\Phi(\bold{u},\boldsymbol\rho)=\sum_{e=1}^n {c}_e(\bold{u})\rho_e
$$
where $c_e=-\rho_e \bold{u}_e^T\bold{K}_e\bold{u}_e<0$. As said in the introduction, we also assume that the elements have the same size. Such assumption is not excessively restrictive and generally used in the literature (see for example \cite{hy10,acs11}). By this assumption we can set $v_e=1,e=1,\dots,n$ and $V=nV_{perc}$, where $V_{perc}$ is the percentage of volume that is occupied by the material in the final design.
With this we can recast problem  (\ref{eq: biknap}) as:

\begin{equation}\label{eq: biles}
\begin{array}{cc}
\displaystyle\max_{\boldsymbol\rho}& \displaystyle\sum_{e=1}^n {c}_e(\bold{u})\rho_e\\
\displaystyle s.t. &\displaystyle \sum_{e=1}^{n} \rho_e\le V\\
&\displaystyle\bold{u}=\mbox{arg} \min_{\bold{v\in{\cal U}_a}}\left\{\frac{1}{2} \bold{v}^T \bold{K}(\boldsymbol\rho)\bold{v}-\bold{v}^T\bold{f}\right\}
\end{array}.
\end{equation}
For the solution of problem (\ref{eq: biles}) we an iterative loop scheme based on the volume reduction strategy generally used in ESO type algorithms. Starting from the initial design composed of only solid elements, at the $k^{th}$ iteration we reduce the maximum percentage of volume that can be occupied by the material by a constant $0<\mu<1$ and have:
\begin{enumerate}
\item Given the design at the previous iteration $\boldsymbol\rho^{k-1}$ solve the lower level problem and compute the displacement variable $\bold{u^k}$;
\item given the displacement fields vector $\bold{u^k}$ use the (\ref{eq: c2}) to compute the sensitivity values associated with the design variables $\boldsymbol\rho$ and solve the knapsack problem to find the new design $\boldsymbol\rho^{k}$.
\end{enumerate}

\section {MATLAB Implementation}\label{mat}
In this section the 51 lines MATLAB code (see Appendix A) is described in detail. The code can be called by the MATLAB prompt by means of the following line:

\bigskip\noindent
\texttt{TOSSE(nelx,nely,volfrac,mu)}

\bigskip\noindent
Where \texttt{nelx} and \texttt{nely} are the number of elements in the horizontal and vertical directions respectively, \texttt{volfrac} is the fraction of volume occupied by the solid elements in the final design and \texttt{mu} is the parameter of volume reduction at every iteration.

The code is obtained after making appropriate modifications to TOP88 , the SIMP code introduced in \cite{acs11}. As a matter of facts, the two codes share the
 first 23 lines  and the finite element analysis is performed in the same manner.  We refer the interested reader to  \cite{acs11} for more details. Differently from TOP88, TOSSE does not need any filter preparation, and immediately starts the main loop at line 30. The iteration is initialized at lines 31-32 where the fraction of feasible volume is immediately reduced according to the parameter \texttt{mu}. The finite element analysis is performed in lines 34-37, in the same fashion as TOP88. The FEM computes the displacement vector $\bold{u}$, that is it solves the lower level problem. As the lower level problem is an unconstrained quadratic optimization problem, the vector $\bold{u}$ is found by solving a system of linear equations through the standard \texttt{mldivide} function in matlab at line 37. Once the value of $\bold{u^k}$ has been found,  the sensitivity analysis is performed to find the values of the vector $\bold{c}_e(\bold{u})$ at line 39, where we add to the current values of $\boldsymbol\rho$ a quantity $E_{min}$ in order to reintroduce elements if they are really rewarding. The value of the objective function is computed at line 40. 


Once  $\bold{c}_e$ is computed the upper level optimization can be executed (lines 42-46). As all the elements have the same size the solution of the generally difficult knapsack problem can be found in sub-polynomial time. As a matter of facts, if all the elements have same size, the global optimal solution can be found by first putting in the knapsack the element with the highest value of $c_e$, then the element with the second to highest value and so on until the sack is full. This can be done in MATLAB by using the \texttt{sort} function on the vector  $\bold{c}$ in  descending order (line 42) and then setting to one the elements corresponding to the first $nV_{k}$ components of the sorted vector (line 44), where  $V_{k}$ is the maximum allowed volume at the current iteration. 
After solving the upper level problem, the results for the iteration are printed at line 48. Finally the algorithm stops when there is no change in the design of the structure.

As the solution of the  knapsack problem can be found in sub-polynomial time, the most expensive operation of the method is the solution of the linear system in order to compute the displacement vector. This operation is quite standard in numerical analysis, and it can be performed efficiently to find a solution for problems with millions of variables in a matter of seconds.

\section{Numerical experience}\label{res}
\begin{figure}
\begin{center}
\subfigure[A half MBB  beam ]{\includegraphics[scale=.35]{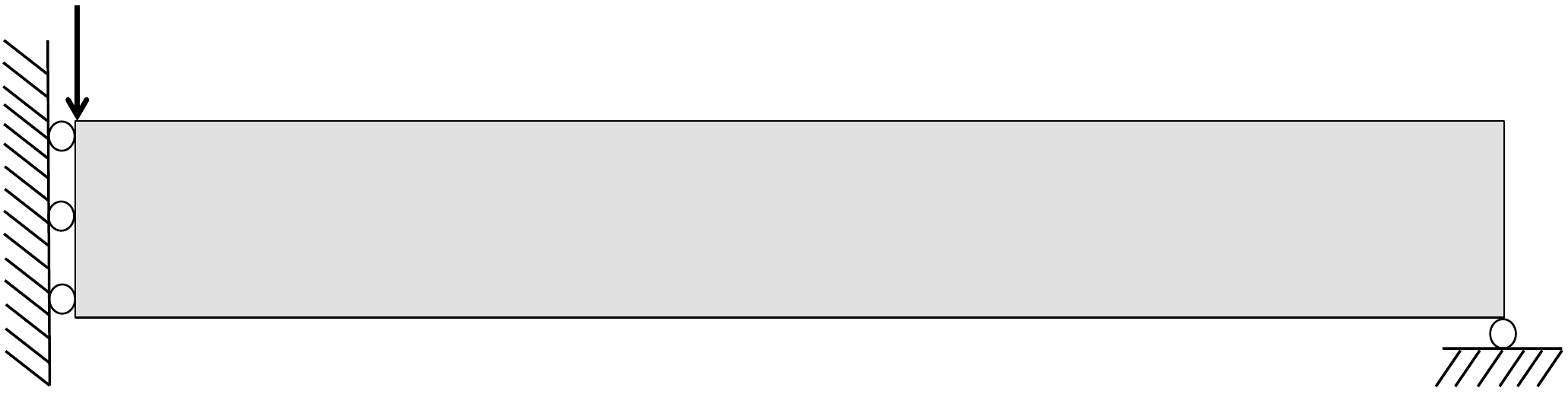}\label{fig:mbb}}
\hspace{.5cm} \subfigure[ A long cantilever   beam ]
{\includegraphics[scale=.35]{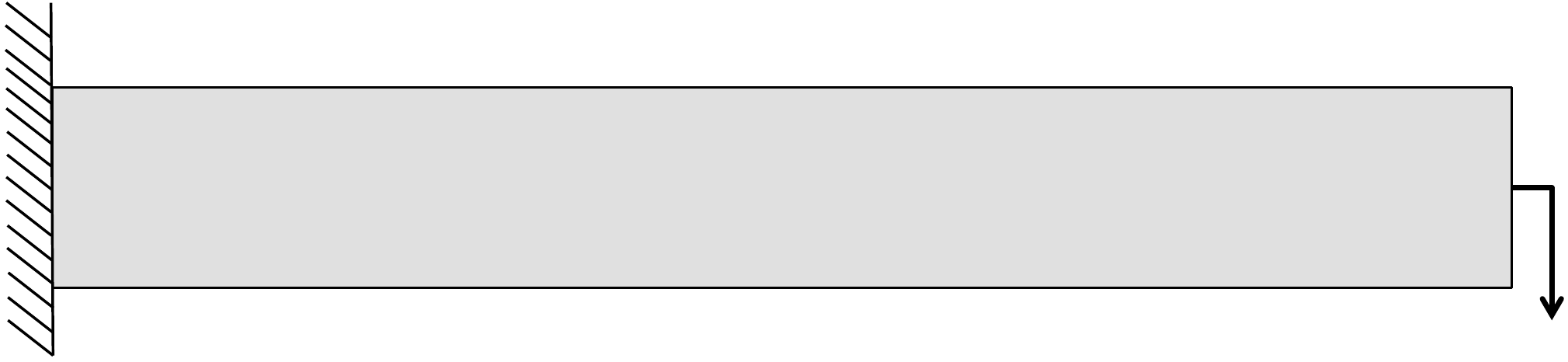}\label{fig:cant}}
\caption{Design domains and external loads for tested beams }\label{fig-beams}
\end{center}
\end{figure}

In this section we report the results of the method proposed so far on the  MBB beam problem in Figure \ref{fig:mbb}.
 These results are intended to show  the practical behavior of the method for general topology problems. In order to better gauge the performances of our algorithm, we also present the numerical results for other codes employed for topology optimization:
\begin{enumerate}
\item TOP88: a popular 88 lines SIMP code for topology optimization code proposed by \cite{acs11};
\item BESO: a soft kill BESO method \cite{hy10};
\item Level Set: an 88 lines  parametrized level set-based topology optimization code proposed in \cite{wei18};
\end{enumerate}
In the first part of this Section we focus on an individual instance of the topology problem to give some qualitative results that are further confirmed in the discussion of the extensive results in Section 4.2. In the results we present we do not use any filter, in order to show how the topologies are  after the optimization and how they compare against checkerboard patters. We remind the reader that checkerboard patters is an undesirable behavior of a topology optimization methods as they are a sign of numerical instability and do not correspond to  optimal distributions of material \cite{check95}.

All the computations in this paper are done using Matlab 2017a on a Linux ubuntu 18.04 64 bits PC with intel(R) core(TM)   I7-6700K CPU and 8 GB of RAM.

\subsection{Qualitative Comparison}
 
\begin{table}[t]
\begin{footnotesize}
\begin{tabular}{c cc>{\centering\arraybackslash} m{6cm}ccc}
\hline
& Optimization  & Total  &Solution & Cost &Time\\
&Parameters   &Iterations&		&Function&\\
\hline
  TOSSE&$\centering \mu=0.97$ &24&\vspace{0.5em}\includegraphics[scale=.4]{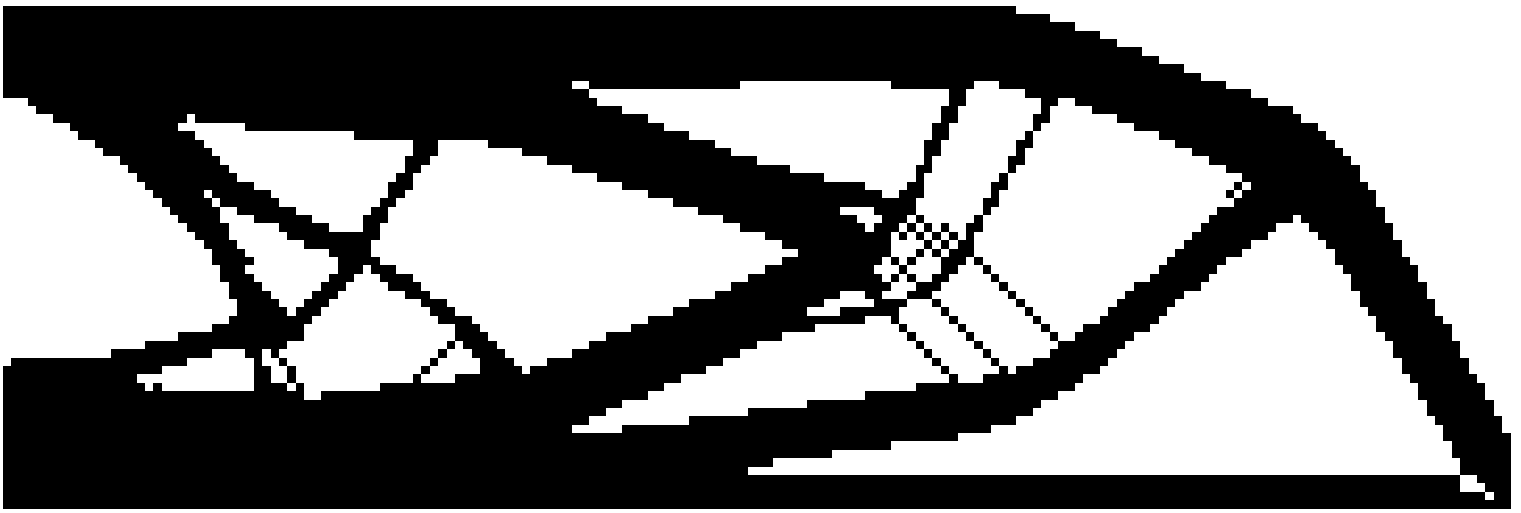}&$C=191.40$&2.3 sec\\
TOP88&\thead{ $p=3$ \\$rmin=1.5$\\ $ft=1$ }&49&\includegraphics[scale=.55]{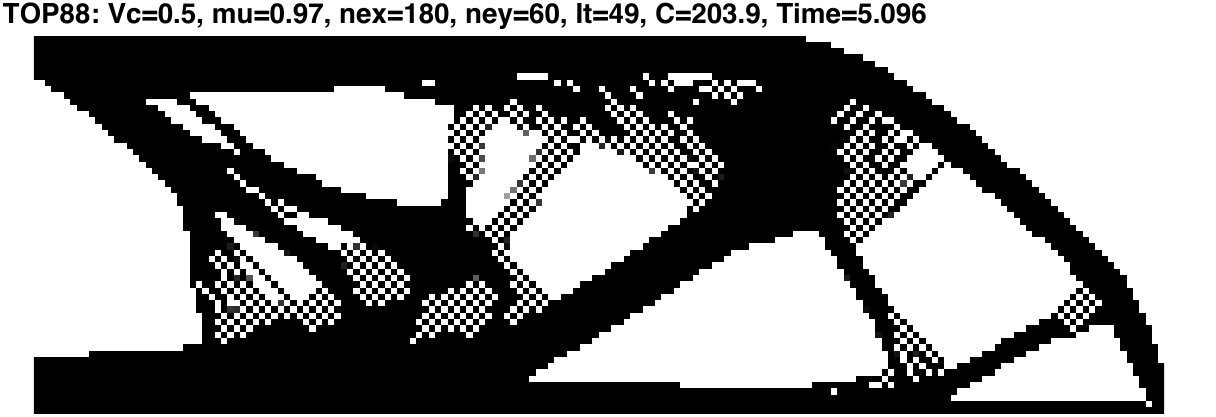}&$C=203.92$&5.1 sec\\
BESO& $\mu=0.97$&28&\includegraphics[scale=.53]{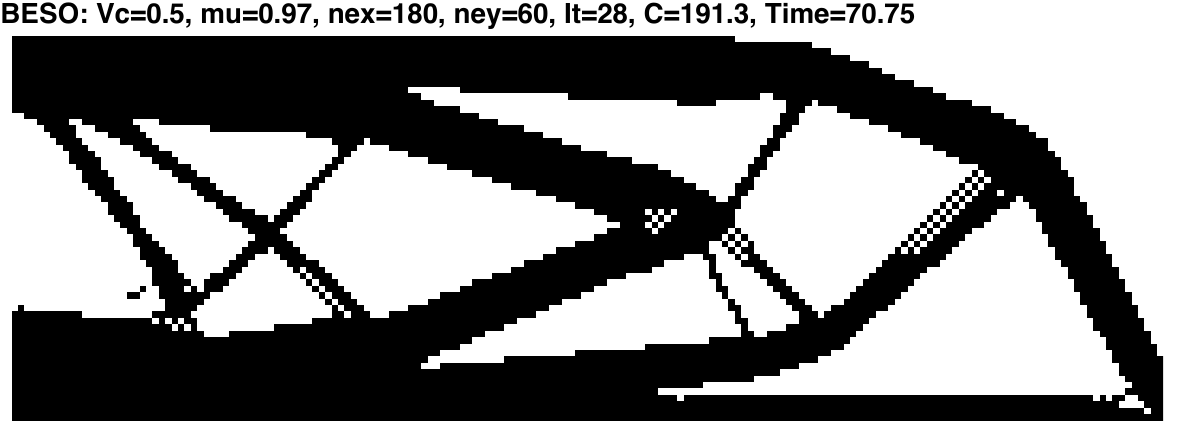}&$C=191.34$&70.7 sec\\
Level Set&\thead{At their\\ default values} &75&\includegraphics[scale=.55]{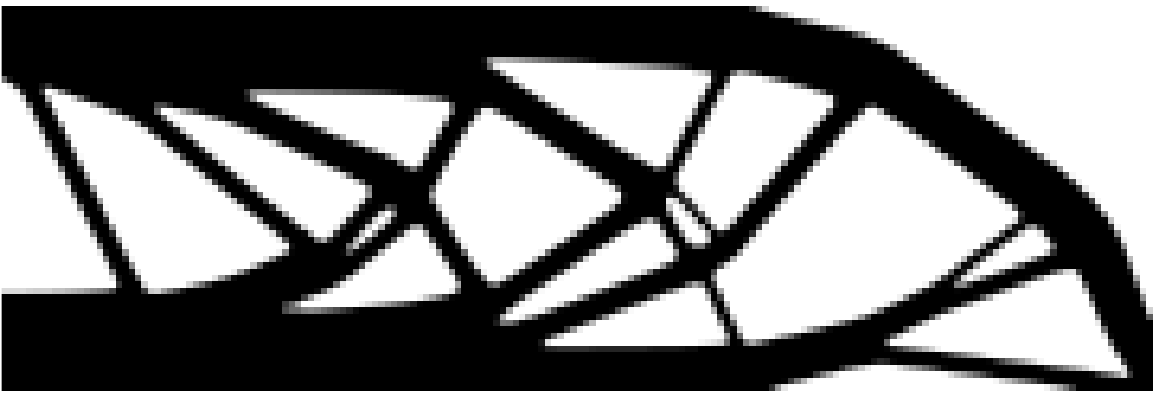}&$C=180.2$&755 sec\\
\hline
\end{tabular}
\end{footnotesize}
\caption{Comparison of the topology optimization methods on the MBB beam.}\label{tab: ex1}
\end{table}

The problem we analyze is a MBB beam with a $180\times60$ design domain. The results are reported in Table \ref{tab: ex1}. It can be easily noticed  that TOSSE and BESO are quite similar in their final designs and values of the objective functions. This should not come to a surprise as both methods use an hard killing evolutionary strategy. The main difference in the design  is that TOSSE shows checkerboard patterns in only one zone of the design while BESO has such problems in four zones. One other important difference is the efficiency. The two methods take almost the same number of iterations to reach the solution, but while the elapsed time per iteration of TOSSE is close to one tenth of second, the same measure is three seconds for BESO. Therefore we can conclude that even if the two designs are quite similar, TOSSE is quite faster than BESO. We will see in the extended results that such difference in efficiency greatly increases with increasing elements.

Regarding efficiency, we see that the one that comes closer to TOSSE in terms of elapsed times is TOP88. The two methods  have quite close elapsed times per iteration, however the final design of TOP88 is characterized by a great number of zones with checkerboard patterns. Therefore TOP88 converges to a local optimum with intermediate elements and with higher values of the objective function. Using a filter would increases the quality of the final solution for TOP88, however the topology obtained by TOSSE could be applied with minor changes in the way it is. Therefore adding a filter would greatly distort the results in favor of TOP88.

Finally the level set method shows a topology without any checkerboard and a low value of the objective function. This is expected, considering the nature of this method. On the other hand this method is the one that takes the highest number of iterations and the most elapsed time to reach a solution, being one order of magnitude slower than BESO that is already an order of magnitude slower than TOSSE and TOP88.

\subsection{Extended Results}
\begin{footnotesize}
\begin{longtable}{|l|l||cc|cc|cc|cc|}
\hline
	&&\multicolumn{2}{c|}{TOSSE}	& \multicolumn{2}{c|}{TOP88}	&\multicolumn{2}{c|}{ BESO}&\multicolumn{2}{c|}{Level Set}	\\
$\mu_c$&nex& Time&Cost&Time&Cost&Time&Cost&Time&Cost\\
\hline
{0.2}
&$180$	&	5.8	&	528.49	&	21.8	&	650.34	&		\multicolumn{2}{c|}{Failure}	&		\multicolumn{2}{c|}{Failure}	\\					
\cline{9-10}																				
&$210$	&	7.7	&	478.33	&	28.4	&	632.17	&	252.8	&	639.42	&		\multicolumn{2}{c|}{\multirow{4}{*}{Out of Memory}}		\\			
&$240$	&	11.5	&	945.49	&	6.1	&	645.27	&		\multicolumn{2}{c|}{Failure}	&		&		\\				
&$270$	&	14.7	&	521.19	&	9.5	&	679.26	&	670	&	515.84	&		&		\\			
&$300$	&	18.2	&	494.16	&	64.4	&	647.18	&	1002.3	&	473.35	&		&		\\			
\hline																				
{0.3}				
&$120$	&	2.0	&	347.25	&	9.1	&	369.64	&		\multicolumn{2}{c|}{Failure}	&		\multicolumn{2}{c|}{Failure}	\\					
&$150$	&	3.0	&	366.66	&	58.4	&	358.03	&	58	&	327.77	&		\multicolumn{2}{c|}{Failure}	\\				
&$180$	&	4.2	&	420.59	&	22.6	&	369.21	&	108.3	&	304.17	&		\multicolumn{2}{c|}{Failure}	\\				
\cline{9-10}																				
&$210$	&	5.9	&	315.99	&	6.4	&	372.64	&	927	&	510.23	&		\multicolumn{2}{c|}{\multirow{4}{*}{Out of Memory}}		\\			
&$240$	&	8.5	&	327.20	&	38.7	&	365.03	&	314.6	&	309.8	&		&		\\			
&$270$	&	12.1	&	367.41	&	20.4	&	372.63	&	528	&	331.24	&		&		\\			
&$300$	&	14.0	&	308.10	&	65.8	&	372.06	&	773	&	313.92	&		&		\\			
\hline																				
{0.5}&$30$	&	0.1	&	195.15	&	0.2	&	190.81	&	0.2	&	234.47	&				\multicolumn{2}{c|}{Failure}	\\		
&$60$ 	&	0.4	&	194.37	&	0.7	&	203.07	&	1.2	&	235.04	&	4.6	&	181.11	\\			
&$90$	&	0.7	&	188.52	&	1.3	&	200.98	&	4.8	&	194.97	&	48.6	&	179.89	\\			
&$120$	&	1.1	&	205.92	&	9.1	&	202	&	14	&	190.77	&	110.9	&	78.61	\\			
&$150$	&	1.7	&	188.16	&	2.7	&	202.74	&	40.1	&	194.41	&	284.9	&	180.5	\\			
&$180$	&	2.3	&	191.40	&	5.1	&	203.92	&	70.7	&	191.34	&	755	&	180.2	\\			
\cline{9-10}																				
&$210$	&	3.5	&	191.63	&	7.2	&	204.33	&	119.6	&	191.63	&		\multicolumn{2}{c|}{\multirow{4}{*}{Out of Memory}}		\\			
&$240$	&	4.6	&	191.81	&	13.9	&	204.13	&	185.6	&	192.01	&		&		\\			
&$270$	&	6.2	&	192.56	&	55.1	&	205.28	&	318.4	&	194.19	&		&		\\			
&$300$	&	7.9	&	193.40	&	66.4	&	204.78	&	484.9	&	197.58	&		&		\\			
\hline																				
{0.7}&$30$	&	0.07	&	143.14	&	0.2	&	143.87	&	0.1	&	143.27	&	0.6	&	142.4	\\			
&$60$ 	&	0.2	&	145.59	&	0.6	&	146.48	&	0.8	&	146.83	&	4.7	&	142.96	\\			
&$90$	&	0.4	&	147.55	&	1.2	&	147.74	&	2.5	&	147.55	&	35.3	&	144	\\			
&$120$	&	0.6	&	151.24	&	1.8	&	148.23	&	7	&	149.77	&	140.2	&	144.93	\\			
&$150$	&	1.0	&	149.77	&	2.8	&	149.05	&	17.3	&	151.59	&	508.8	&	145.78	\\			
&$180$	&	1.3	&	151.62	&	6.9	&	149.98	&	38	&	153.52	&	1564.4	&	146.34	\\			
\cline{9-10}																				
&$210$	&	1.9	&	151.75	&	7.5	&	150.36	&	59.6	&	153.13	&		\multicolumn{2}{c|}{\multirow{4}{*}{Out of Memory}}		\\			
&$240$	&	2.5	&	152.03	&	12.7	&	151	&	108.1	&	153.82	&		&		\\			
&$270$	&	3.4	&	152.54	&	13.5	&	151.41	&	184.3	&	153.35	&		&		\\			
&$300$	&	4.9	&	153.01	&	64.9	&	151.77	&	238.8	&	154.04	&		&		\\			
\hline
\caption{Results for the four algorithms on the MBB problem, increasing the number of elements and the volume fraction. The value of ney is one third of nex.}\label{tab: mbb}
\end {longtable}
\end{footnotesize}

In this section we present the results on an extended variety of tests. 
Through the numerical experience we change two parameters:
\begin{itemize}
\item The number of elements on the mesh $nex$ and $ney$; 
\item The volume fraction $\mu_c$. This parameter determines the maximum fraction of volume that can be filled by the elements. 
\end{itemize}
Increasing the number of elements on the mesh brings to more complex and robust structures, but also increases the computational cost, as the number discrete variables goes with the square of the number of elements on the axes. With many  elements we can see how  the algorithms behave on large sized optimization problems. Furthermore the lower the volume fraction is, more the optimization problem is difficult as some algorithms could fail to create a reasonable structure.
In our testing the number of elements on the y axis is one third of the number of elements on the x axis, and we increase the number of elements on the x axis of 30 units (10 on the y axis) in every optimization until we reach 300 elements. For $\mu_c=0.2$ we start with $nex=180$ elements, for $\mu_c=0.3$ we start with $nex=120$ elements, while for $\mu_c=0.5,0.7$ the initial number of elements is $nex=30$.

The results are reported in Table \ref{tab: mbb} and they confirm the qualitative considerations made in the previous subsection. As a matter of facts it is clear that the  fastest method on average is TOSSE, with TOP88 having comparable elapsed times, while BESO is tens of times slower and the Level Set method is even slower. It is interesting to underline that TOSSE and TOP88 have quite close elapsed times per iteration, however TOP88 takes more iterations to reach a solution. We also notice that the relative growth of the elapsed times for TOSSE is the smallest among all the algorithms.

For what regards the cost function, once again TOSSE and BESO reach quite similar results especially in topologies with many elements while TOP88 generally has higher values of the cost and the Level Set method lower values. The only exceptions for TOP88 are the topologies with $\mu_c=0.7$, however the designs returned by such algorithm in these cases are characterized by a large amount of checkerboard patterns.

\begin{table}[t]
\begin{tiny}
\begin{tabular}{>{\centering\arraybackslash}m{4.4cm}>{\centering\arraybackslash}m{3.5cm}>{\centering\arraybackslash}m{3.5cm}>{\centering\arraybackslash}m{3.5cm}}
$\mu_c$=0.2\includegraphics[scale=.25]{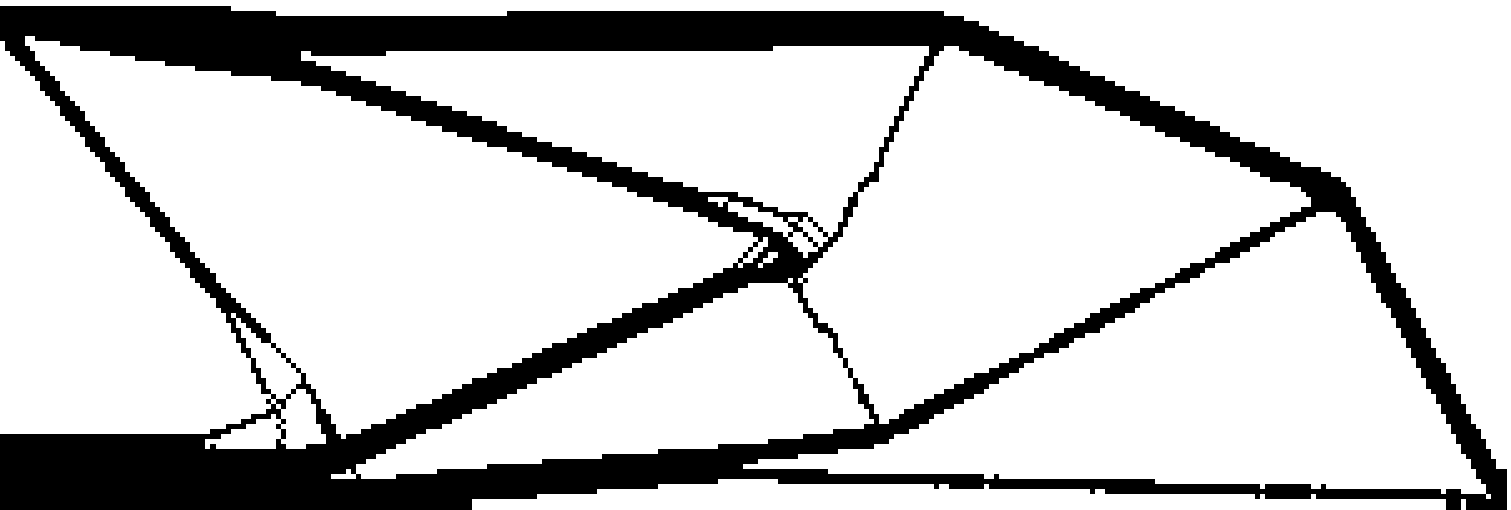}&\includegraphics[scale=.33]{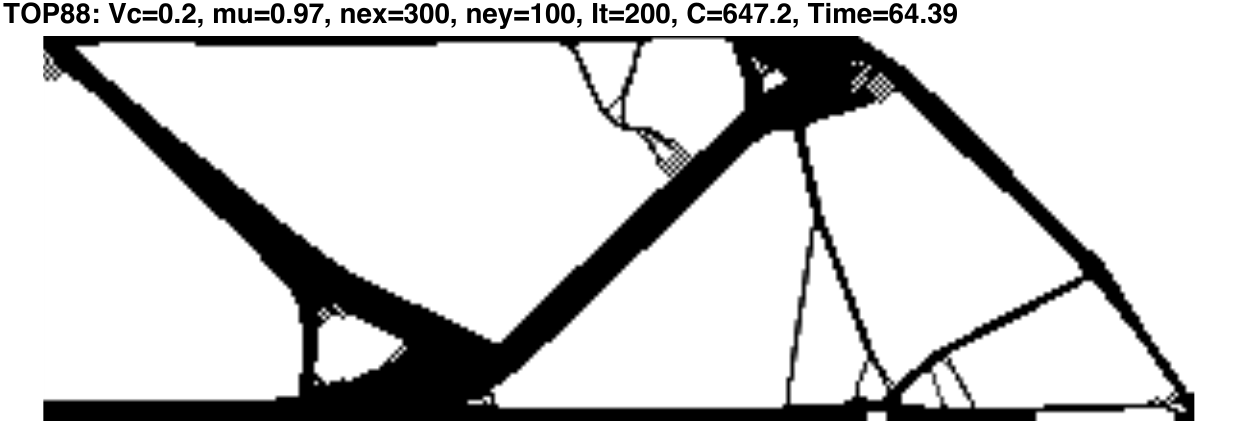}&\includegraphics[scale=.33]{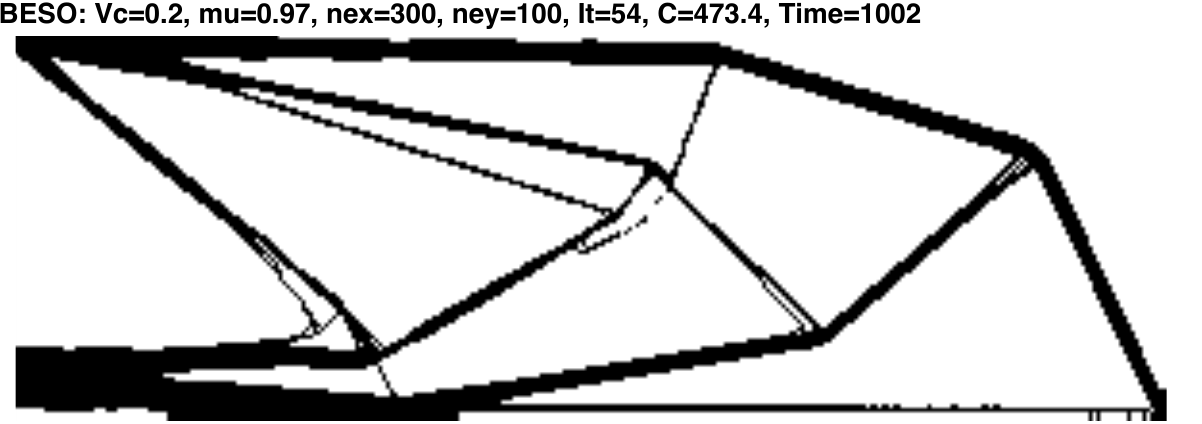}&\\
$\mu_c$=0.3\includegraphics[scale=.25]{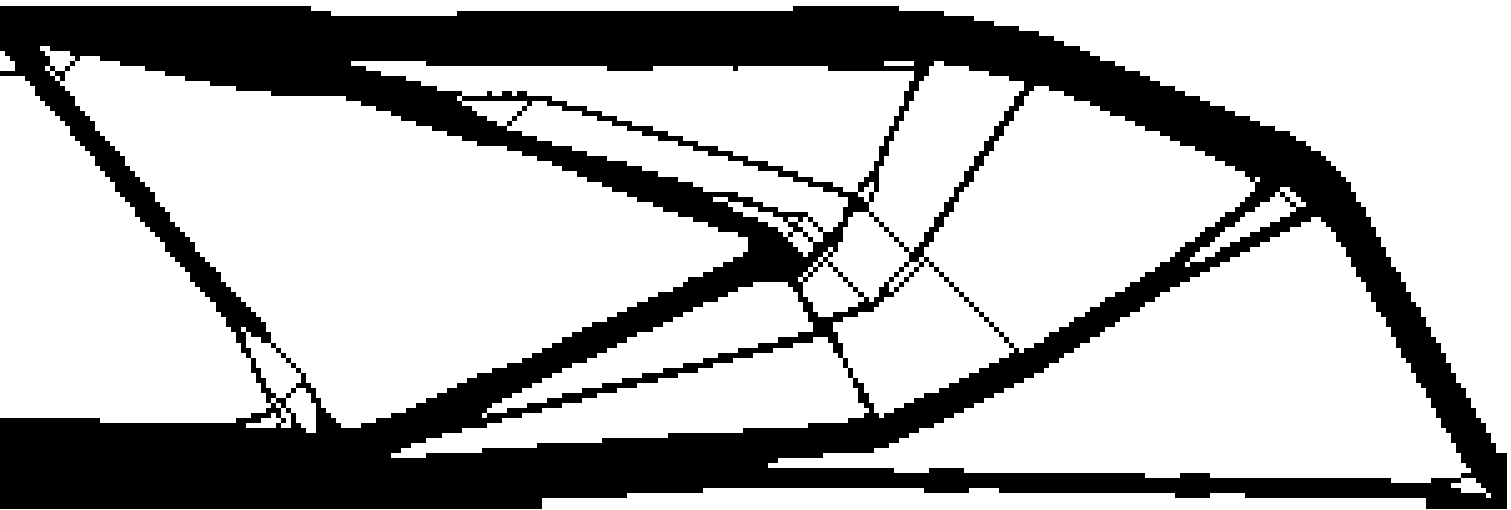}&\includegraphics[scale=.33]{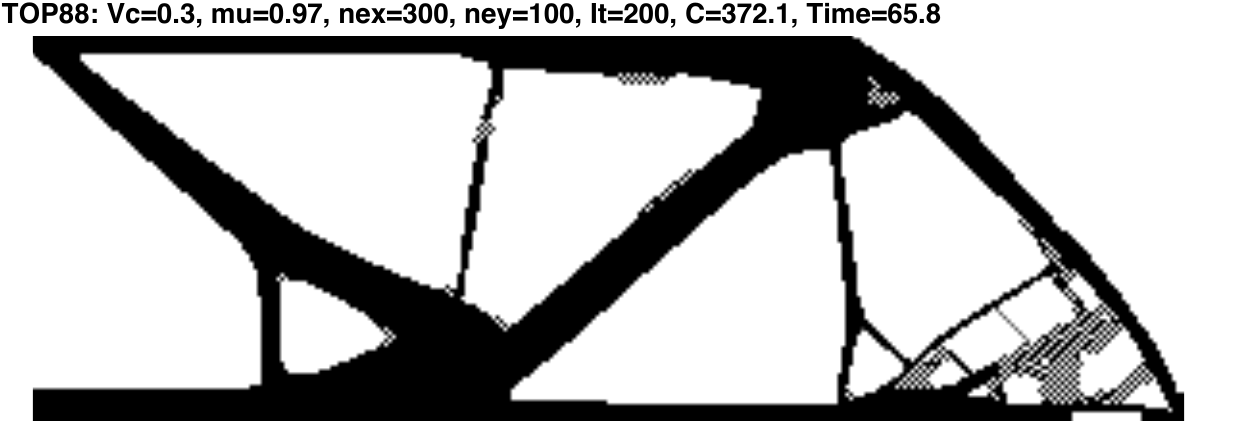}&\includegraphics[scale=.33]{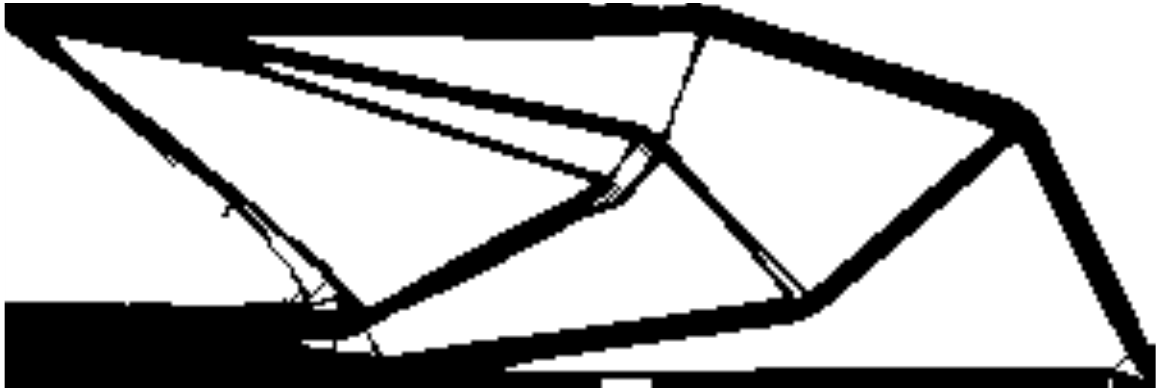}&\\
$\mu_c$=0.5\includegraphics[scale=.33]{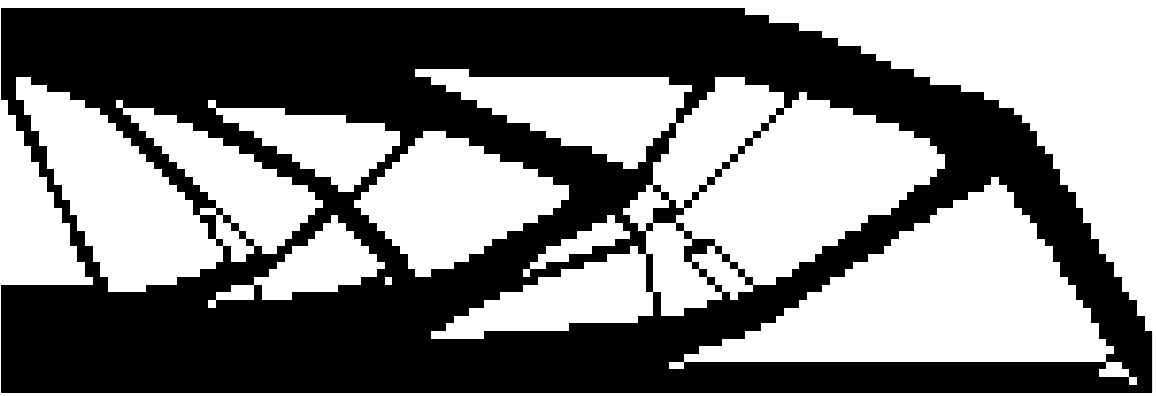}&\includegraphics[scale=.33]{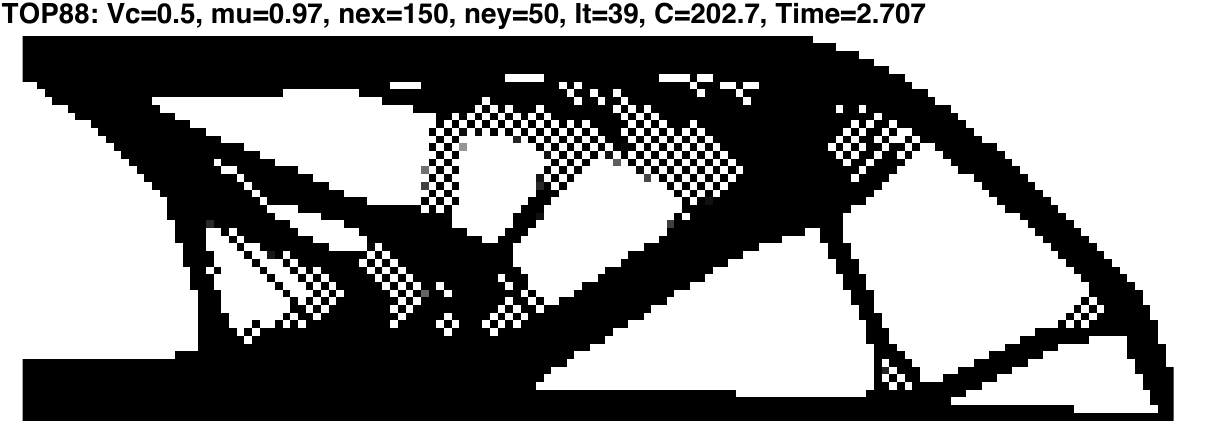}&\includegraphics[scale=.33]{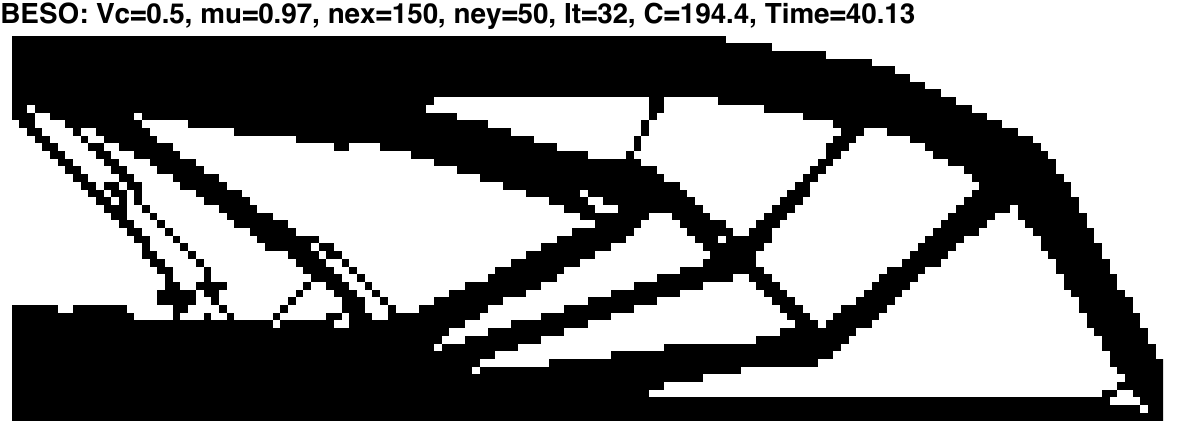}&\includegraphics[scale=.33]{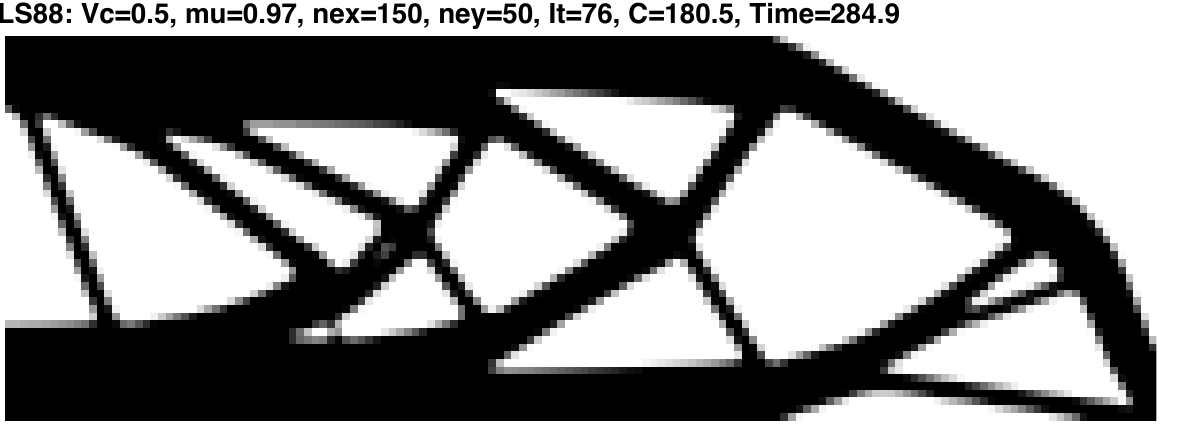}\\
$\mu_c$=0.7\includegraphics[scale=.25]{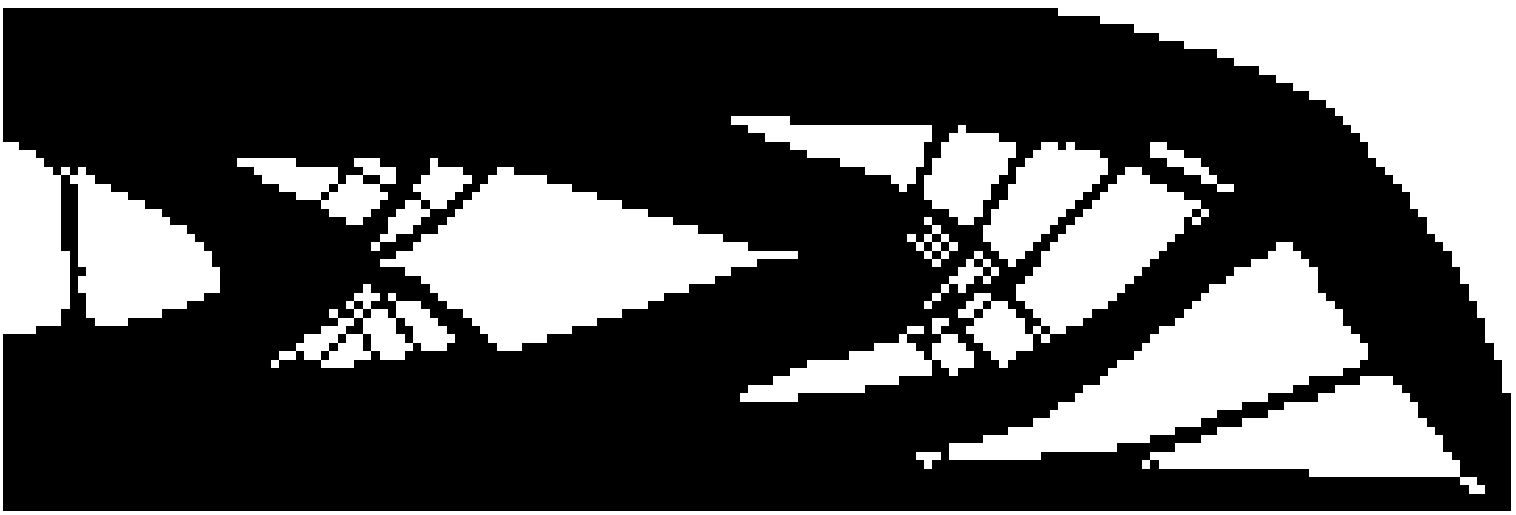}&\includegraphics[scale=.33]{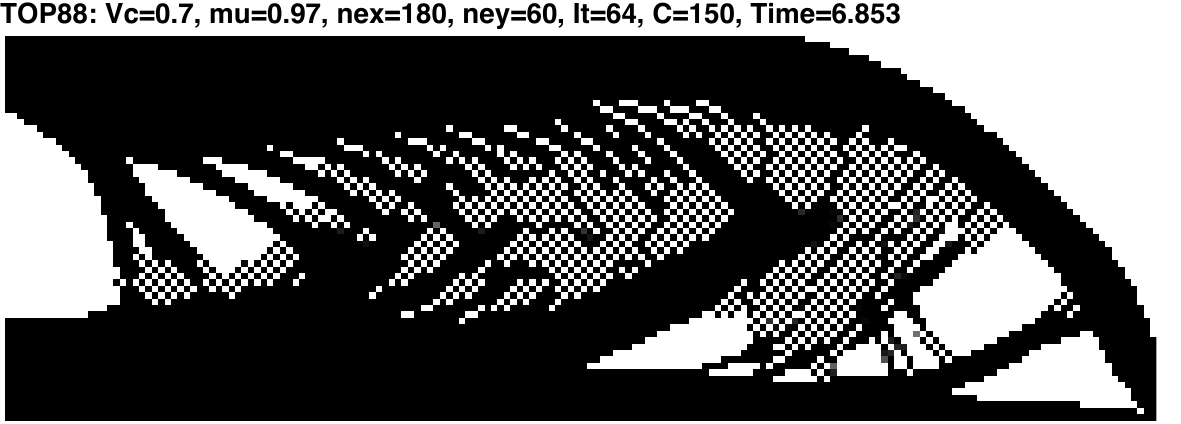}&\includegraphics[scale=.33]{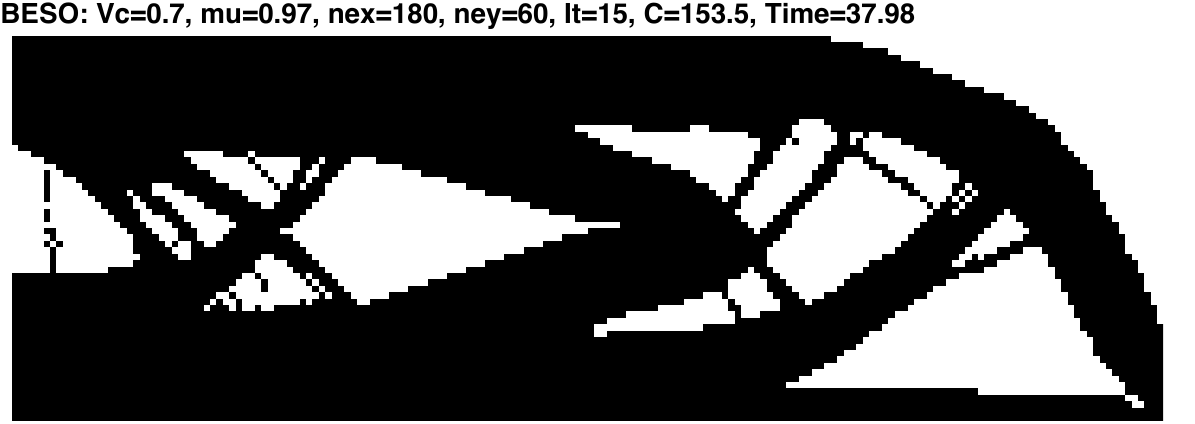}&\includegraphics[scale=.33]{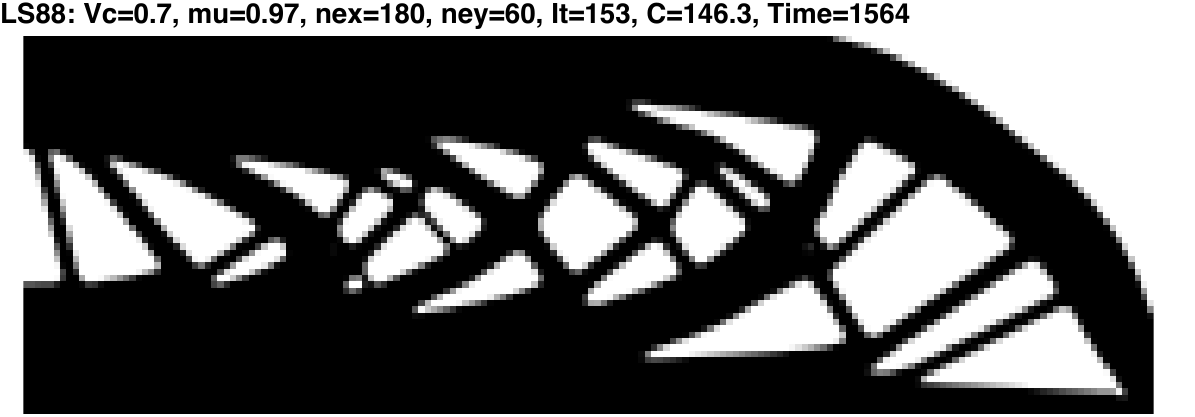}\\
TOSSE&TOP88&BESO&Level Set\\
\end{tabular}
\end{tiny}
\caption{MBB beam with  $\mu_c=0.2,0.3,0.5, 0.7$.}\label{tab: mbb_fig}
\end{table}

In Table \ref{tab: mbb_fig} we report the topologies for the different algorithms and different values of the $\mu_c$ parameters. For $\mu_c=0.2,0.3$ the Level Set method is not able to yield any topology, while for $\mu_c=0.5,0.7$ the maximum number of elements that the Level Set method is able to handle is  $nex=180$ and $ney=60$. Therefore we present design with $nex=300$ and $ney=100$ for $\mu_c=0.2,0.3$, $nex=150$ and $ney=50$ for $\mu_c=0.5$ (The final designs for $nex=180$ and $ney=60$ are already reported in Table \ref{tab: ex1}) and $nex=180$ and $ney=60$  for $\mu_c=0.7$.
These final topologies confirm the consideration made in Section 4.1. TOSSE and BESO generate topologies that are quite similar and with almost no checkerboard patterns, differently by TOP88 that is characterized by many zones with checkerboards. The structures yielded by TOSSE are also simpler than those obtained by the Level Set method and characterized by less connections.

Finally we make some considerations on the stability of the different softwares and the memory usage. For what regards the stability, both TOSSE and TOP88 behave in similar manners and are more stable than BESO. The level set method instead is unable to generate structures with small values of  $\mu_c$ and small values of $nex$ and $ney$. It probably would be able to generate a topology with designs with an high number of elements, but such method requires too much memory to complete this task. BESO  requires more memory than TOSSE and TOP88, while these last two methods have more or less the same memory usage.
As a matter of the facts the most memory taxing operation for the two codes is the creation of the sparse matrix to be used to compute the displacements at every iteration. The solution of the linear system to compute the displacements  is also the most time expensive task of both the algorithms and this explains why TOSSE and TOP88 have similar elapsed times per iteration.

To summarize the results of this analysis, we can say that TOSSE is the fastest method on average, capable of yielding topologies with almost no checkerboard patters. It is on average faster than TOP88 and yields topologies similar to BESO. Furthermore both the reliability and the memory usage of the software are on the level of TOP88 showing that it can solve difficult topology optimization problems.

\section{Extensions}\label{ext}

TOSSE can be extended in several manners, and in this Section we show how to use it on the two dimensional and three dimensional cantilever beams. We also offer the codes that implement these variations on the webpage \url{https://github.com/vlatorre847/TOSSE}{}. We underline that many of the extensions proposed in \cite{sig99,acs11} can be easily applied to this code with minor changes.

\subsection{2D Cantilever Beam}
Another well known problem in topology optimization is the cantilever beam (Figure  \ref{fig:cant}). In order to solve this problem in topology optimization,  the following line must be added after line 20:

\bigskip\noindent
\texttt{F(2*((nely+1)*nelx+(ceil(nely/2)+1)),1) = -1; }

\bigskip\noindent
And the boundary conditions at line 21 must be changed to

\bigskip\noindent
\texttt{fixeddofs=[1:2*(nely+1)];}

\bigskip\noindent
Comparisons has been made with the other algorithms similarly to the results described in Section 4. In Table \ref{tab: can} we report the numerical comparisons with the other topology algorithms for $\mu_c=0.5$.
\begin{table}
\begin{footnotesize}
\begin{center}
\begin{tabular}{|l|l||cc|cc|cc|cc|}
\hline
	&&\multicolumn{2}{c|}{TOHE}	& \multicolumn{2}{c|}{TOP88}	&\multicolumn{2}{c|}{ BESO}&\multicolumn{2}{c|}{Level Set}	\\
$\mu_c$&nex& Time&Cost&Time&Cost&Time&Cost&Time&Cost\\\hline
0.5&$30$	&	0.1	&	427.50	&	1.4	&	188.12	&		\multicolumn{2}{c|}{Failure}	&	5.3	&	164.81	\\
&$60$ 	&	0.3	&	179.67	&	2.6	&	186	&	1.2	&	180.58	&	21.6	&	164.49	\\
&$90$	&	0.6	&	187.93	&	3.5	&	182.82	&	4.7	&	222.9	&	79.3	&	164.63	\\
&$120$	&	1.1	&	177.22	&	2.3	&	183.41	&	13.2	&	183.32	&	325.1	&	165.17	\\
&$150$	&	1.9	&	177.04	&	8.4	&	182.82	&	30.3	&	179.27	&	807.8	&	165.35	\\
\cline{9-10}																	
&$180$	&	2.4	&	183.95	&	6.3	&	184.36	&	60.5	&	173.04	&		\multicolumn{2}{c|}{\multirow{5}{*}{Out of Memory}}	\\	
&$210$	&	3.4	&	180.79	&	30.1	&	184.65	&	110.6	&	174.31	&		&		\\
&$240$	&	4.9	&	173.16	&	40.2	&	184.59	&	193.5	&	176.48	&		&		\\
&$270$	&	7.1	&	176.61	&	55.2	&	184.7	&	305.4	&	174.54	&		&		\\
&$300$	&	8.8	&	171.45	&	17.2	&	185.24	&	445.6	&	174.27	&		&		\\
\hline																	
\end{tabular}
\end{center}
\caption{Results for the four algorithms on the Cantilever problem, increasing the number of elements for volume fraction $\mu_c=0.5$. The value of ney is one third of nex.}\label{tab: can}
\end{footnotesize}
\end {table}
The results in Table \ref{tab: can} are a further confirmation of the considerations made in the previous Section, that is TOSSE is the fastest method on average with cost function similar to BESO and with almost no checkerboard patterns on the final topology. 

One of the things that make TOSSE differ from the other methods is that the final topologies yielded by this method for the cantilever beam could  be not symmetric, that is the upper part of the beam differs from the lower part. This is due to the sensitivities $c_e$, that also represent the cost in the knapsack problem, are not perfectly symmetric. This could be due to numerical errors in the solution of the linear system that computes the displacements. It is possible to avoid this issue by solving the knapsack problem in a standard fashion and then reflect the upper or lower half of the topology in order to obtain a symmetric design. 
\begin{figure*}
\begin{tabularx}{0.5\textwidth}{cc}
\hspace{-1.3em}\subfigure[No symmetry]{
\begin{tabular}{c}
\small $\mu_c$=0.2{\includegraphics[scale=.45]{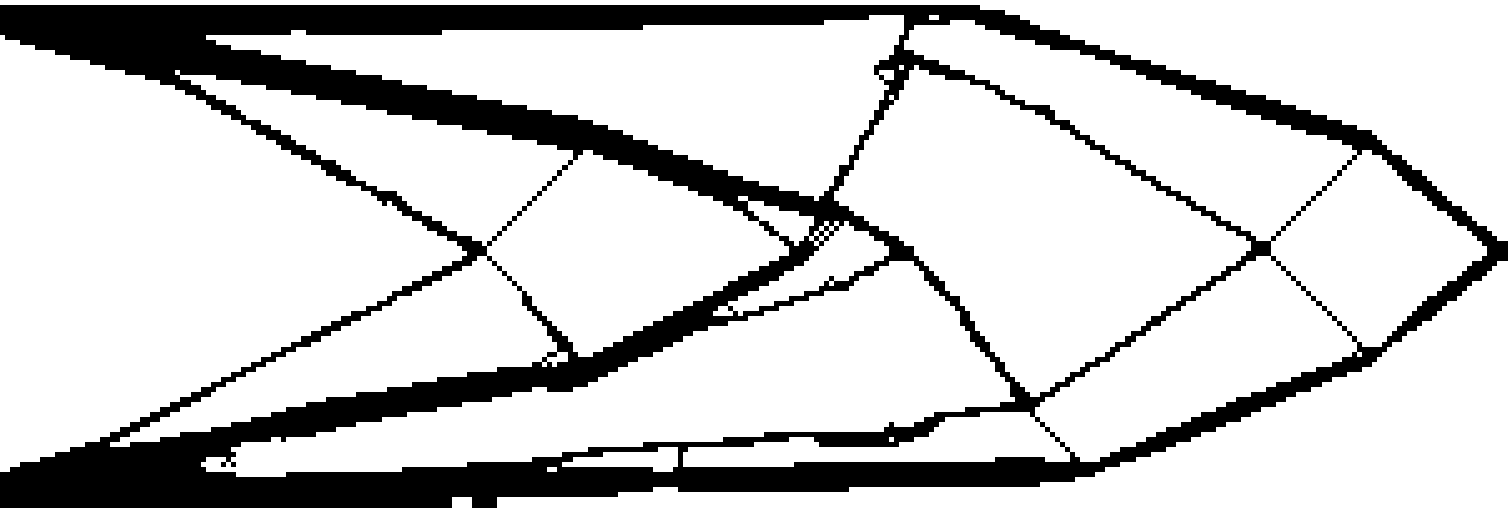}}\\
\small $\mu_c$=0.3{\includegraphics[scale=.45]{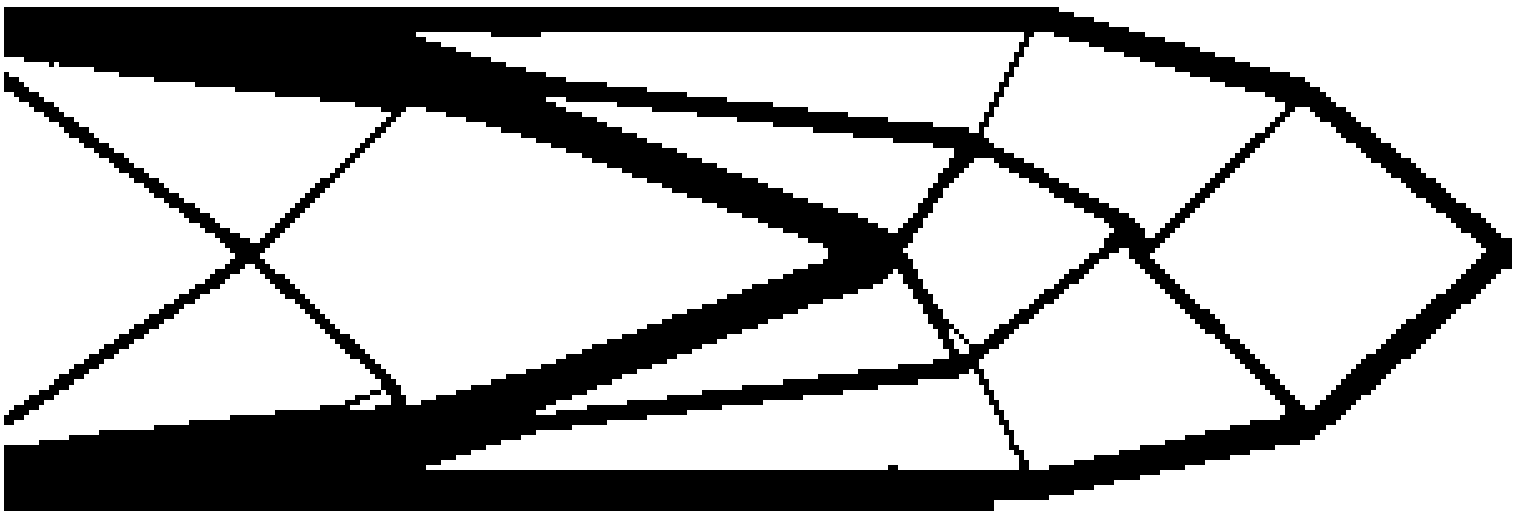}}\\
\end{tabular}}&
\hspace{-2.5em}\subfigure[Symmetry]{
\begin{tabular}{c}
{\includegraphics[scale=.45]{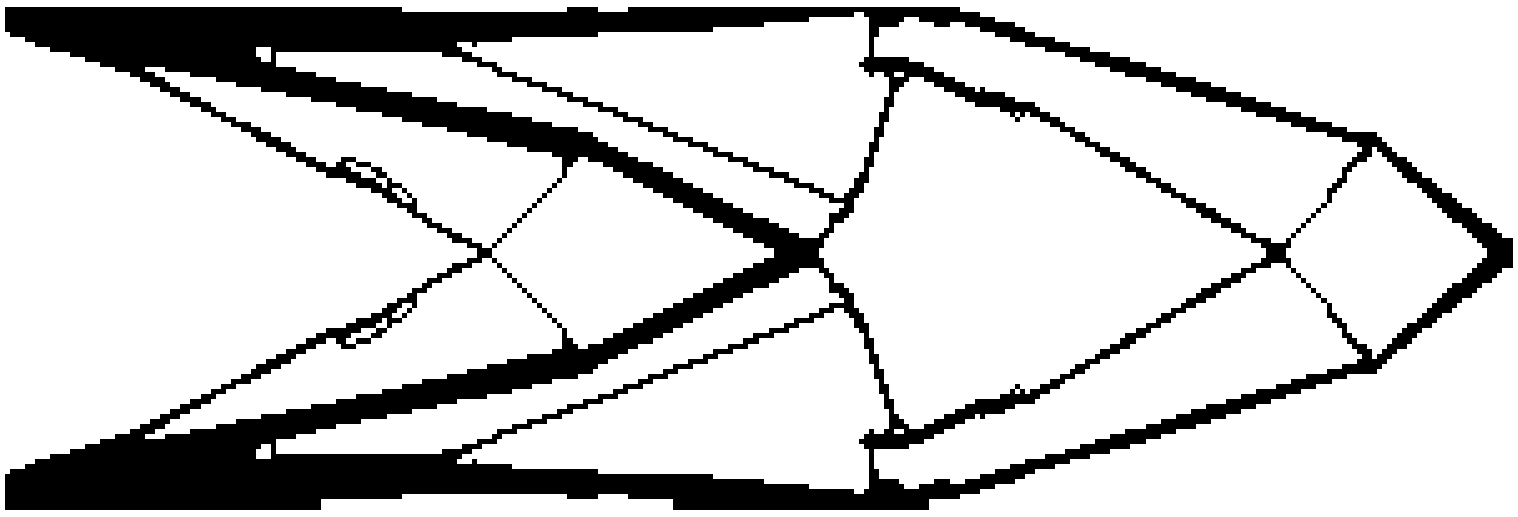}}\\
{\includegraphics[scale=.45]{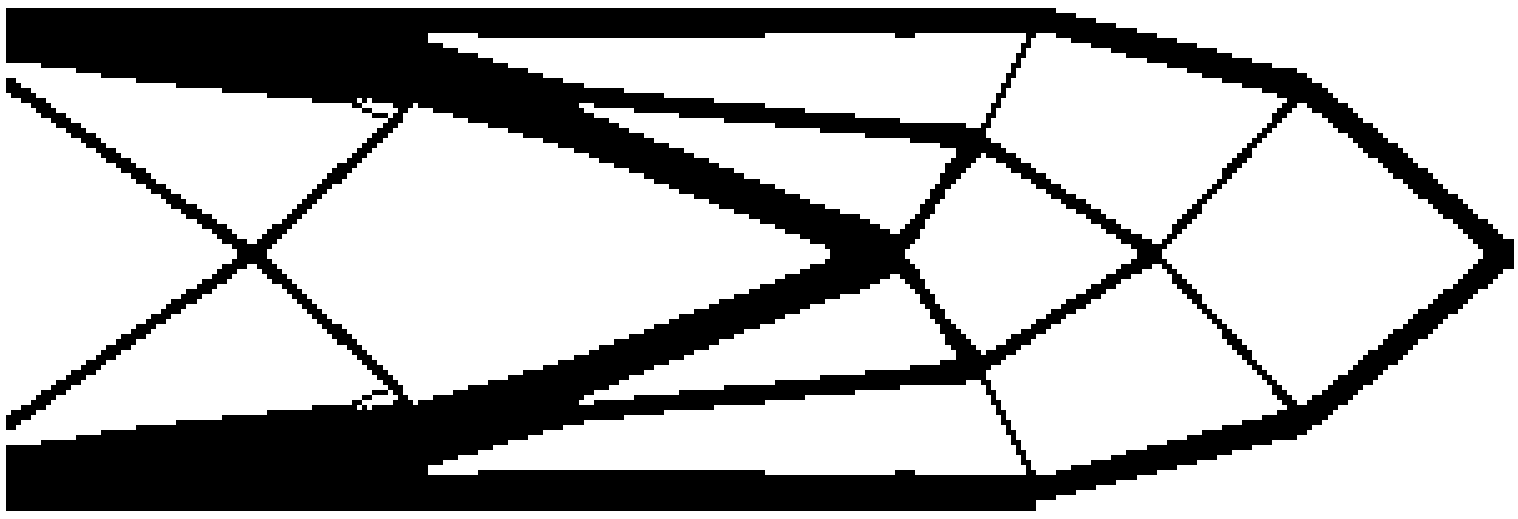}}\\
\end{tabular}}
\end{tabularx}
\caption{The behavior of TOSSE for the cantilever beam with symmetrical and asymmetrical strategies, $nelx=300$, $nely=100$. }\label{fig: ss}
\end{figure*}
In Figure \ref{fig: ss} it is possible to see the comparison between the  asymmetric topologies and the symmetric topologies.  The  topologies do not differ  much, however there is the possibility that asymmetric topologies are better, as the asymmetry in the sensitivities could be due to an asymmetrical distribution of the forces on the single elements rather than to numerical errors.
\subsection{Three dimensional structures}
\begin{figure*}
\begin{tabularx}{0.5\textwidth}{cc}
\hspace{-1.3em}\subfigure[TOSSE3D]{
\begin{tabular}{c}
\small{\includegraphics[scale=.45]{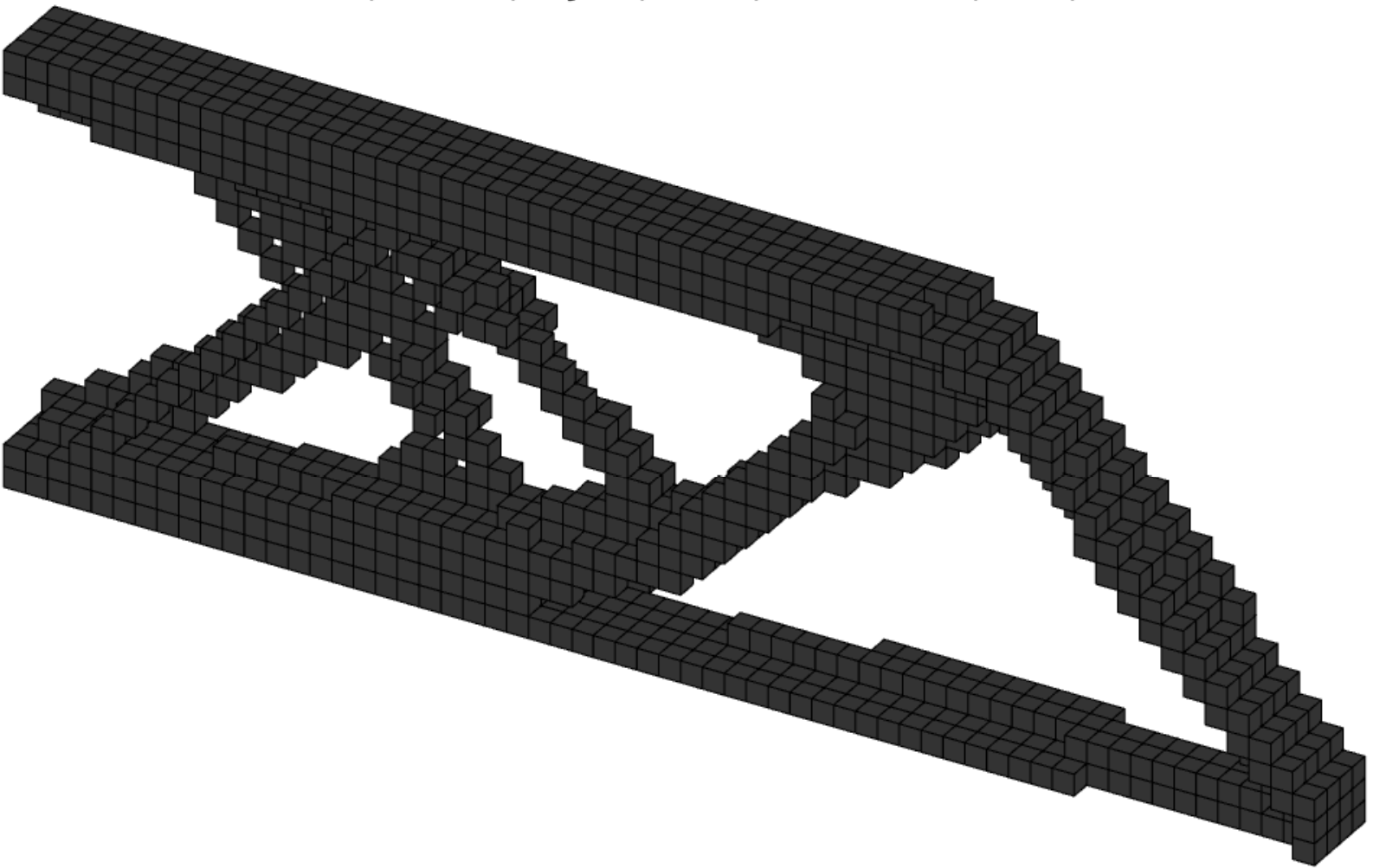}}
\end{tabular}}&
\hspace{-2.5em}\subfigure[TOP3D]{
\begin{tabular}{c}
{\includegraphics[scale=.45]{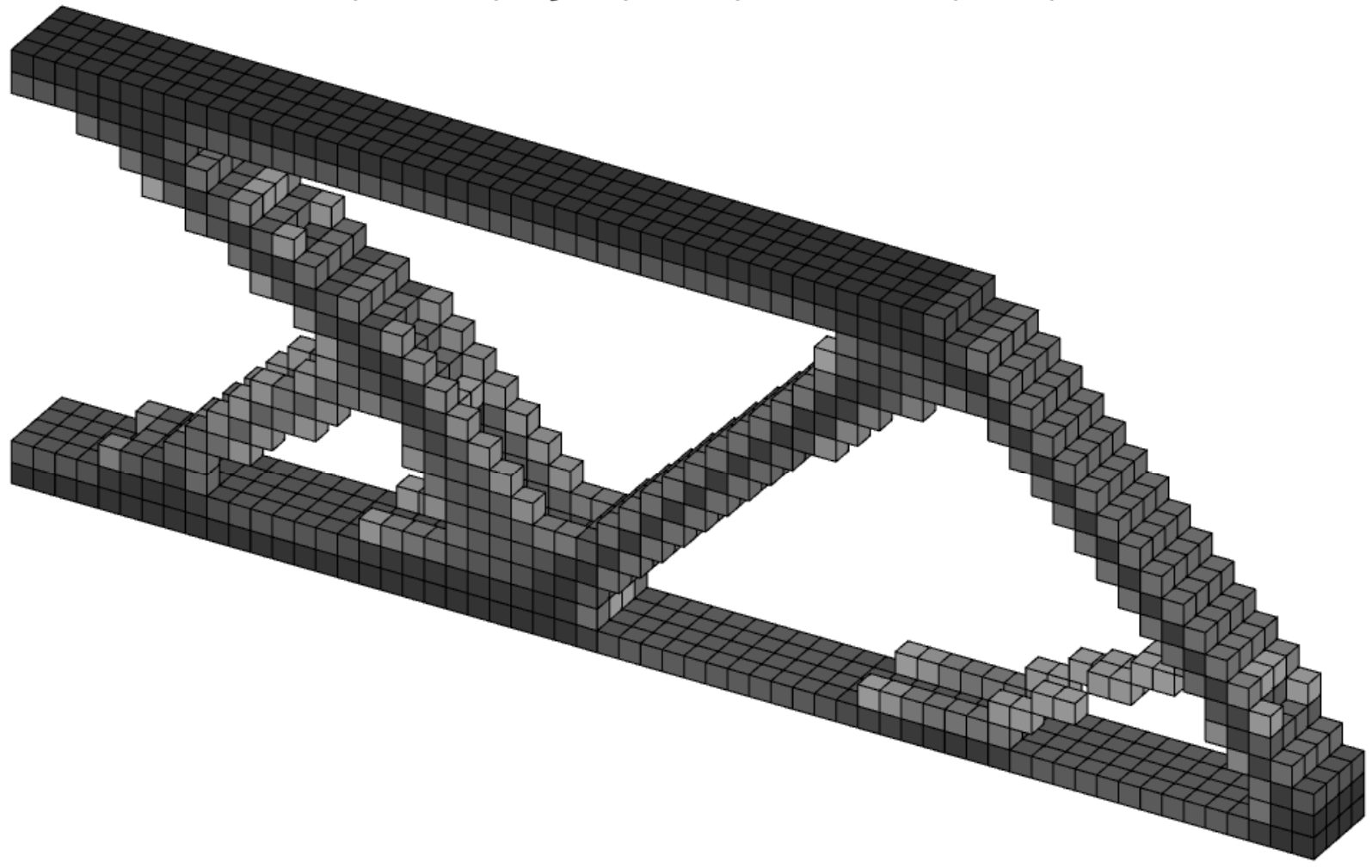}}\label{top3d}
\end{tabular}}
\end{tabularx}
\caption{Comparisons between TOSSE3D and TOP3D for three dimensional structures. The lighter the grey color is, the less the element is close to the integer 1 value.}\label{fig: 3d}
\end{figure*}

The application of the knapsack strategy to three dimensional problems is straightforward. The code we propose is based on the code reported in \cite{lt14} that is a 169 lines Matlab code for three dimensional topology optimization problems that exploit a SIMP strategy. The changes we made to the code are quite similar to the two dimensional case. In Figure \ref{fig: 3d} we report the results for the 3D cantilever beam on a topology with a $60\times20\times4$ design domain and $\mu_c=0.3$. It is possible to see that in the three dimensional case TOP3D has a clear difficulty to put elements to either 0 or 1, as showed in the Figure \ref{top3d}, where the  lighter shades of gray of many of the elements   indicate that their values are far from 1. This three dimensional version of TOP does not have checkerboard patterns, but instead is not able to get clear void-solid design that cannot be applied without the use of a filter. 
On the other hand, TOSSE3D creates a topology with clear void-solid elements that satisfies the volume constraints and can be directly applied with small modifications. 
From the point of view of elapsed times, the two method have once again the same elapsed times per iteration, with TOSSE3D being able to get a solution in quite less iterations.

\section{Conclusions}\label{con}
This paper presents a new Matlab code for the topology optimization problem. This new methodology is based on the formulation of the problem as a bilevel mixed integer optimization problem and the assumption that all elements have the same size.  

The code is 51 lines long, can be considered a very simple implementation of a method for the solution of the topology optimization problem, and can be used as a basis for further developments of methods in this field. The new code has been written in a way that does not compromise its readability and efficiency although being quite shorter than the other codes present in the literature. Therefore it is suitable for educational purposes.

From the numerical experience, it emerges that the efficiency of this code is on average faster than the classical SIMP method, with almost none checkerboard patterns and clear void-solid structures. 
The results indicate that this new method seems to yield final designs similar to those generated by the evolutionary methods such as ESO/BESO  but with the efficiency of SIMP methods.

The Matlab code can be easily extended to several other boundary conditions and even to three dimensional problems. We argue that the many extensions already existing in the literature can be easily adapted to this code. On the Github page \url{https://github.com/vlatorre847/TOSSE}{} associated with this paper, it is possible to download the code for the MBB and cantilever beams for 2D structures as well as the code for 3D cantilever structures.
This code, based on a novel formulation and efficient implementation, can be the basis for new developments and extensions in the field of topology optimization.

\section*{Acknowledgments}
The research was  supported by  US Air Force Office of  Scientific Research
 under the grant  FA9550-17-1-0151.

\appendix
\section{Matlab Code}
\lstset{language=Matlab,%
    breaklines=true,%
    morekeywords={matlab2tikz},
    keywordstyle=\color{blue},%
    morekeywords=[2]{1}, keywordstyle=[2]{\color{black}},
    identifierstyle=\color{black},%
    stringstyle=\color{mylilas},
    commentstyle=\color{mygreen},%
    showstringspaces=false,
    numbers=left,%
    numberstyle={\tiny \color{black}},
    numbersep=9pt, 
}
\begin{lstlisting}
%%% A 51 line code for topology optimization with same size elements. Vittorio Latorre Dec. 2018 %%%
function [c,xnew,loop]=tosse(nelx,nely,volfrac,mu)
%% MATERIAL PROPERTIES
E0 = 1;
Emin = 1e-9;
nu = 0.3;
%% PREPARE FINITE ELEMENT ANALYSIS
A11 = [12  3 -6 -3;  3 12  3  0; -6  3 12 -3; -3  0 -3 12];
A12 = [-6 -3  0  3; -3 -6 -3 -6;  0 -3 -6  3;  3 -6  3 -6];
B11 = [-4  3 -2  9;  3 -4 -9  4; -2 -9 -4 -3;  9  4 -3 -4];
B12 = [ 2 -3  4 -9; -3  2  9 -2;  4  9  2  3; -9 -2  3  2];
KE = 1/(1-nu^2)/24*([A11 A12;A12' A11]+nu*[B11 B12;B12' B11]);
nodenrs = reshape(1:(1+nelx)*(1+nely),1+nely,1+nelx);
edofVec = reshape(2*nodenrs(1:end-1,1:end-1)+1,nelx*nely,1);
edofMat = repmat(edofVec,1,8)+repmat([0 1 2*nely+[2 3 0 1] -2 -1],nelx*nely,1);
iK = reshape(kron(edofMat,ones(8,1))',64*nelx*nely,1);
jK = reshape(kron(edofMat,ones(1,8))',64*nelx*nely,1);
% DEFINE LOADS AND SUPPORTS (HALF MBB-BEAM)
U = zeros(2*(nely+1)*(nelx+1),1);
F = sparse(2,1,-1,2*(nely+1)*(nelx+1),1);
fixeddofs = union([1:2:2*(nely+1)],[2*(nelx+1)*(nely+1)]);
alldofs = [1:2*(nely+1)*(nelx+1)];
freedofs = setdiff(alldofs,fixeddofs);
%% INITIALIZE ITERATION
vgam = 1;
loop = 0;
change = 1;
xPhys = ones(nely,nelx)*vgam;
%% START ITERATION
while change >   1e-16 &&loop<200
    loop = loop + 1;
    vgam = max(volfrac,vgam*mu);
    %% FE-ANALYSIS
    sK = reshape(KE(:)*(Emin+xPhys(:)'.*(E0-Emin)),64*nelx*nely,1);
    K = sparse(iK,jK,sK);
    K = (K+K')/2;
    U(freedofs) = K(freedofs,freedofs)\F(freedofs);
    %% OBJECTIVE FUNCTION AND SENSITIVITY ANALYSIS
    ce = (xPhys(:)+Emin)*E0.*(sum((U(edofMat)*KE).*U(edofMat),2));
    c = sum((Emin+ce*(1-Emin/E0)));
    %% SOLVE THE KNAPSACK PROBLEM
    [~,I]=sort(ce,1,'descend');
    xnew=zeros(nelx*nely,1);
    xnew(I(1:floor(vgam*nelx*nely)))=1;
    change = max(abs(xnew(:)-xPhys(:)));
    xPhys = reshape(xnew,[nely,nelx]);
    %% PRINT RESULTS
    fprintf(' It.:%5i Obj.:%11.4f Vol.:%7.3f ch.:%7.3f fea.:%7.3f\n',loop,c,mean(xPhys(:)),change,sum(sum(xPhys))-vgam*nelx*nely);
end
%% PLOT DENSITIES
colormap(gray); imagesc(1-xPhys); caxis([0 1]); axis equal; axis off; drawnow;
\end{lstlisting}

\bibliographystyle{plain}

\end{document}